\documentclass[12pt]{amsart}
\usepackage{amssymb,amsmath,amsfonts,latexsym}
\usepackage{amsfonts,amssymb,amscd,amsmath,enumerate,verbatim,calc}
\usepackage{bm}

\newcommand{\op}{\displaystyle\oplus_{1}^n}
\newcommand{\q}{\bigoplus_{i=1}^n}

\newcommand{\sm}{\sum_{j=1}^n}

\newcommand{\ov}{\overline}

\newcommand{\z}{\mathcal B}
\newcommand{\p}{^\prime}
\newcommand{\m}{\mathcal}

\theoremstyle{plain}
\newtheorem{theorem}{Theorem}[section]

\newtheorem{lemma}[theorem]{Lemma}

\numberwithin{equation}{section}

\begin{document}

\title{Factorizations of Characteristic Functions}

\author[Haria]{ Kalpesh J. Haria}
\address{Indian Statistical Institute, Statistics and Mathematics Unit, 8th Mile, Mysore Road, Bangalore, 560059, India}
\email{kalpesh\_if@isibang.ac.in, hikalpesh.haria@gmail.com}

\author[Maji]{Amit Maji}
\address{Indian Statistical Institute, Statistics and Mathematics Unit, 8th Mile, Mysore Road, Bangalore, 560059, India}
\email{amaji\_pd@isibang.ac.in, amit.iitm07@gmail.com}

\author[Sarkar]{Jaydeb Sarkar}
\address{Indian Statistical Institute, Statistics and Mathematics Unit, 8th Mile, Mysore Road, Bangalore, 560059, India}
\email{jay@isibang.ac.in, jaydeb@gmail.com}



\begin{abstract}
Let $A = (A_1, \ldots, A_n)$ and $B = (B_1, \ldots, B_n)$ be row
contractions on $\mathcal{H}_1$ and $\mathcal{H}_2$, respectively,
and $X$ be a row operator from $\oplus_{i=1}^n \mathcal{H}_2$ to
$\mathcal{H}_1$. Let $D_{A^*} = (I - A A^*)^{\frac{1}{2}}$ and
$D_{B} = (I - B^* B)^{\frac{1}{2}}$ and $\Theta_T$ be the
characteristic function of $T =
\begin{bmatrix}
A& D_{A^*}L D_B\\
0 & B \end{bmatrix}$. Then $\Theta_T$ coincides with the product of
the characteristic function $\Theta_A$ of $A$, the Julia-Halmos
matrix corresponding to $L$ and the characteristic function
$\Theta_B$ of $B$. More precisely, $\Theta_T$ coincides with
\[
\begin{bmatrix}
\Theta_B & 0 \\
0 & I
\end{bmatrix}
(I_\Gamma \otimes \begin{bmatrix} L^* & (I - L^* L)^{\frac{1}{2}}
\\ (I - L L^*)^{\frac{1}{2}} & - L
\end{bmatrix})
\begin{bmatrix}
\Theta_A & 0\\ 0& I\end{bmatrix},
\]
where $\Gamma$ is the full Fock space. Similar results hold for
constrained row contractions.
\end{abstract}

\maketitle

\section{Introduction}

Let $\m H$ be a Hilbert space and $T$ be a contraction (that is, $I
- T T^* \geq 0$) on $\m H$ and suppose $D_T = (I - T^*
T)^{\frac{1}{2}}$ and $D_{T^*} = (I - T T^*)^{\frac{1}{2}}$ are the
defect operators and $\m D_T = \overline{\mbox{ran}} D_T$ and $\m
D_{T^*} = \overline{\mbox{ran}} D_{T^*}$ are the defect spaces of
$T$. Then the \textit{characteristic function} of $T$ is an operator
valued bounded analytic function $\Theta_T \in H^\infty_{\m B(\m
D_{T}, \m D_{T^*})}(\mathbb{D})$ defined by
\[
\Theta_T(z) = (-T + D_{T^*} z (I - z T^*)^{-1} D_T)|_{\m D_T} \quad
\quad (z \in \mathbb{D}).
\]
The notion of characteristic functions plays an important role in
many areas of operator theory and function theory (see
\cite{NFBK10}). In particular, characteristic functions are one of
the central objects of study in noncommutative operator theory and
noncommutative function theory (see Popescu \cite{P10} and
references therein).

On the other hand, the notion of Julia-Halmos matrix is important in
the construction of isometric and unitary dilation maps for
contractions (cf. \cite{NFBK10}). Recall that the
\textit{Julia-Halmos matrix} corresponding to a contraction $L$ from
$ \m H $ to $ \m K$ is the unitary matrix
\[
J_L = \begin{bmatrix} L^* & (I - L^* L)^{\frac{1}{2}}
\\ (I - L L^*)^{\frac{1}{2}} & - L
\end{bmatrix}
=
\begin{bmatrix} L^*& D_{ L}\\
D_{ L^*} &- L
\end{bmatrix}.
\]
This is also directly related to analytic or functional models for
contractions in the sense of Sz.-Nagy and Foias (see Timotin
\cite{Tim}). However, one of the most striking results along these
lines is due to Sz.-Nagy and Foias (see \cite{NF67}):

\noindent\textsf{Theorem (Sz.-Nagy and Foias):} Let $\m H_1$ and $\m
H_2$ be Hilbert spaces and $T =
\begin{bmatrix}T_1 & X\\ 0 & T_2\end{bmatrix}$ be a contraction on $\m H_1 \oplus
\m H_2$. Then there exist a contraction $L \in \m B(\m D_{T_1^*},
\m D_{T_2})$ and (canonical) unitary operators $\tau \in \m B(\m D_T, \m
D_{T_1} \oplus \m D_{L})$ and $\tau_* \in \m B(\m D_{T^*}, \m
D_{T_2^*} \oplus \m D_{L^*})$ such that $X = D_{T_1^*} L D_{T_2}$
and
\[\Theta_T(z) = \tau_*^{-1} \begin{bmatrix} \Theta_{T_2}(z) & 0 \\ 0
& I_{\m D_{L^*}}\end{bmatrix}
\begin{bmatrix} L^* & D_{L}
\\ D_{L^*} & - L
\end{bmatrix}
\begin{bmatrix} \Theta_{T_1}(z)
& 0 \\ 0 & I_{\m D_L}\end{bmatrix} \tau \quad \quad (z \in
\mathbb{D}).\]

In this paper we first generalize the above factorization result to
noncommuting tuples of row contractions. For the class of
constrained row contractions, we obtain a similar result to the main
factorization result.

The paper is organized as follows: In Section 2 we give a brief
introduction of characteristic functions and multi-analytic
functions in noncommutative set up and fix some notations. In
Section 3 we present the Sz.-Nagy and Foias type factorization
results for noncommuting tuples of row contractions. In the final
section we obtain a similar factorization result for constrained row
contractions.

\section{Preparatory results}

In this section we recall and study some basic tools of
multivariable operator theory such as characterizations of upper
triangular operator matrices, characteristic functions and
multi-analytic functions which appear in all later investigation. A
general theory of characteristic operators and (multi-)analytic
models for row contractions on Hilbert spaces was developed by G.
Popescu in \cite{Po89b}, \cite{Po89c} and \cite{Po99} (also see
\cite{Po06a} and references therein).

Let $\ H$ be a Hilbert space and $\{T_j\}_{j=1}^n \subseteq \z (\m
H)$ where $\z (\m H)$ is the algebra of all bounded linear operators
on $ \m H $. Then the $n$-tuple $T = (T_1, \ldots, T_n) $ is called
a row contraction if $T : \q \m H \rightarrow \m H$ is a
contraction, that is, $\sm T_j T_j^* \leq I_{\m H}$ or,
equivalently,  if $\|\sm T_j h_j\|^2 \leq \sm \|h_j\|^2$, $h_1,
\ldots, h_n \in \m H$.

The defect operators and defect spaces of a row contraction $T$ on
$\m H$ are given by
\[D_{T} = (I- T^* T)^\frac{1}{2} \in \z(\q\m H), \quad \quad D_{{T}^*} = (I- TT^*)^\frac{1}{2} \in \z(\m
H),
\]
and
\[\m D_{T} = \overline{\mbox{ran}}D_T \subseteq \q\m H \quad \quad \m D_{T^*} = \overline{\mbox{ran}}
D_{T^*} \subseteq \m H,
\]
respectively.

The class of row contractions with which we concerned has the
following characterization (see \cite{NF67} or Lemma 2.1, Chapter IV
in \cite{FF90}):

\begin{theorem}\label{th:NFblock}\textsf{(Sz.-Nagy and Foias)}
Let $\m H_1$ and $\m H_2$ be two Hilbert spaces and $A =  (A_1,
\ldots, A_n)\in \z (\op \m H_1, \m H_1)$, $B = (B_1, \ldots, B_n)
\in \z( \op \m H_2, \m H_2)$ and $X = (X_1, \ldots, X_n) \in \z(\op
\m H_2,\m H_1 )$ are row operators. Then the row operator
\[T =
\begin{bmatrix}A & X\\ 0 & B\end{bmatrix} \in \m B( (\op \m H_1) \oplus(\op \m H_2), \m H_1 \oplus \m
H_2)
\]
is a row contraction if and only if $A$ and $B$ are row contractions
and
\[X = D_{A^*} L D_{B},
\]
for some contraction $L \in  \z (\m D_{B}, \m D_{A^*})$.
\end{theorem}

We now recall the following result of Sz.-Nagy and Foias about
unitary operators between defect spaces (see \cite{NF67} or
Corollary 2.2, Chapter IV in \cite{FF90}):

\begin{theorem}\label{th:NFdefect}\textsf{(Sz.-Nagy and Foias)}
In the setting of Theorem \ref{th:NFblock}, let $T$ be a row
contraction. Then there exist unitary operators  $\sigma: \m D_{T}
\to \m D_{A} \oplus \m D_L $ and $ \sigma_*: \m D_{T^*} \to \m
D_{B^*} \oplus \m D_{L^*} $ such that
\begin{align}\label{2}
\sigma D_T = \begin{bmatrix}
D_{A}&  - A^* L D_{B}\\
0 & D_L D_{B}\\
\end{bmatrix} \quad \mbox{and} \quad  \sigma_* D_{T^*}= \begin{bmatrix}
- B L^*D_{A^*}&  D_{B^*}\\
D_{L^*}D_{A^*} & 0\\
\end{bmatrix}.
\end{align}
\end{theorem}

The full \textit{Fock space} over $\mathbb C^n$, denoted by
$\Gamma$, is the Hilbert space
\[
\Gamma
:=\displaystyle\bigoplus_{m=0}^{\infty}(\mathbb{C}^n)^{\otimes^m}=
\mathbb{C}\oplus\mathbb{C}^n\oplus(\mathbb{C}^n)^{\otimes^2} \oplus
\cdots\oplus(\mathbb{C}^n)^{\otimes^m}\oplus \cdots.
\]
The \textit{vacuum vector} $1\oplus 0 \oplus \cdots \in \Gamma$ is
denoted by $e_\emptyset$. Let $\{e_1,\ldots, e_n\}$ be the standard
orthonormal basis of $\mathbb C^n$ and $\mathbb{F}_n^+$ be the
unital free semi-group with generators $1,\ldots, n$ and the
identity $\emptyset$. For $\alpha = \alpha_1\ldots\alpha_m \in
{\mathbb{F}_n^+}$ we denote the vector $e_{\alpha_1}\otimes \ldots
\otimes e_{\alpha_m} \in \Gamma$ by $e_\alpha$. Then  $\{e_\alpha:
\alpha \in {\mathbb{F}_n^+}\}$ forms an orthonormal basis of
$\Gamma$. For each $j = 1, \ldots, n$, the left creation operator
$L_j$ and the right creation operator $R_j$ on $\Gamma$ are defined
by
\[L_j f = e_j \otimes f, \quad \quad R_j f = f \otimes e_j \quad \quad (f \in
\Gamma),
\]
respectively. Moreover, $ R_j = U^* L_j U $ where $U$, defined by
$U(e_{i_1}\otimes e_{i_2}\otimes\ldots \otimes e_{i_m}) =  e_{i_m}
\otimes \ldots \otimes e_{i_2}\otimes e_{i_1} $, is the flipping
operator on $\Gamma$. The \textit{noncommutative disc algebra} ${\m
A_n^\infty}$ is the norm closed algebra generated  by $\{I_{\Gamma},
L_1, \ldots,L_n \}$ and the \textit{noncommutative analytic Toeplitz
algebra} $\m F^\infty_n $ is the WOT-closure of ${\m A_n^\infty}$
(see Popescu \cite{Po95a}).

Let $\m E$ and $\m E_*$ be Hilbert spaces and $M \in \m B(\Gamma
\otimes \m E,  \Gamma \otimes \m E_*)$. Then $M$ is said to be
\textit {multi-analytic operator} if
\[
 M(L_j \otimes I_{\m E} ) = (L_j \otimes I_{\m E_*}) M  \quad \quad (j =1, \ldots,n).
\]
In this case the bounded linear map $\theta \in \m B(\m E, \Gamma
\otimes \m E_*)$ defined by
\[\theta (\eta) = M(e_\emptyset \otimes \eta) \quad \quad (\eta \in
\m E),
\]
is said to be the \textit{symbol} of $M$ and we denote $M =
M_\theta$. Moreover, define $\theta_{\alpha} \in \m B(\m E, \m
E_*)$, $\alpha \in {\mathbb{F}_n^+} $ by
\[
\langle  \theta_{\alpha} \eta, \eta_* \rangle := \langle \theta
\eta, e_{\bar \alpha} \otimes \eta_* \rangle = \langle M
(e_\emptyset\otimes \eta),  e_{\bar \alpha} \otimes \eta_* \rangle,
\quad \quad (\eta \in \m E, \eta_* \in \m E_*)
\]
where $\bar \alpha$ is the reverse of $ \alpha$. Then one can associate a unique formal
power series
\[
M ~\sim ~ \displaystyle\sum_{\alpha \in {\mathbb{F}_n^+}} R_\alpha
\otimes \theta_{\alpha},
\]
and (see Popescu \cite{Po95a})
\[
M = \mbox{SOT}-\displaystyle\lim_{r \to 1^{-}}
\displaystyle\sum_{k=0}^{\infty}\displaystyle\sum_{|\alpha|=k}r^{|\alpha|}
R_\alpha \otimes \theta_{\alpha}
\]
where $ |\alpha|$ is
the length of $ \alpha$.

A multi-analytic operator $M_{\theta} \in \m B (\Gamma \otimes \m E,
\Gamma \otimes \m E_*)$  is said to be \textit{purely contractive}
if $M_{\theta}$ is a contraction and
\[
\|P_{ e_\emptyset \otimes \m E_*} \theta \eta \| < \|\eta\| \quad
\quad (\eta \in \m E, \eta \neq 0).
\]

We say that $M_\theta$ coincides with a multi-analytic operator
$M_{\theta\p} \in \m B(\Gamma  \otimes \m E\p, \Gamma  \otimes \m
E_*\p)$ if there exist unitary operators $ W: \m E \to \m E\p $ and
$ W_*:\m E_* \to \m E_*\p  $ such that
\[
(I_{\Gamma}\otimes W_*) M_{\theta} =  M_{\theta\p}
(I_{\Gamma}\otimes W).
\]

Let $\m H$ be a Hilbert space and  $T = (T_1, \ldots, T_n)$ be a row
operator on $\m H$. \textsf{For simplicity of the notations we will
denote by $\tilde T$ and $\tilde R$ the row operators $(I_\Gamma
\otimes T_1, \ldots, I_{\Gamma} \otimes T_n)$ and $(R_1 \otimes
I_{\m H}, \ldots, R_n \otimes I_{\m H})$ on $\Gamma \otimes \m H$,
respectively.}

Among multi-analytic operators, characteristic functions
\cite{Po06b} play an important role in multivariable operator
theory and noncommutative function theory (see \cite{P10} and other
references therein). The \textit{characteristic function} of a row
contraction $T$ on $\m H$ is a purely contractive multi-analytic
operator $\Theta_{T} \in \m B(\Gamma \otimes \m D_{T}, \Gamma
\otimes \m D_{T^*})$ defined by
\[
\Theta_T \sim~ -I_\Gamma \otimes T + (I_\Gamma \otimes  D_{
T^*})(I_{\Gamma \otimes \m H}- \tilde R \tilde T^*)^{-1} \tilde R
(I_\Gamma\otimes D_T).
\]
Hence
\[
\Theta_T  = \mbox{SOT}-\displaystyle\lim_{r \to 1}\Theta_T(r R),
\]
where for each $r \in [0, 1)$,
\[
\Theta_T(r R) := - \tilde T + D_{ {\tilde T}^*}(I_{\Gamma \otimes \m
H}-  r \tilde R\tilde T^*)^{-1} r \tilde R  D_{\tilde T}.
\]
Therefore
\begin{align}\label{3}
\Theta_T = \mbox{SOT}-\displaystyle\lim_{r \to 1}\Theta_T(r R) =
\mbox{SOT}-\displaystyle\lim_{r \to 1} \big[- \tilde T+ D_{ {\tilde
T}^*}(I_{\Gamma \otimes \m H}-  r \tilde R\tilde T^*)^{-1} r \tilde
R  D_{\tilde T}\big].
\end{align}

\section{Factorizations of Characteristic Functions of noncommuting tuples}

In this section we prove the main theorem on factorizations of
characteristic functions of upper triangular operator matrices. We
begin with the following simple lemma.

\begin{lemma}\label{le:theta}
Let $T$ be a row contraction on $\m H$. Then for each $r \in [0, 1)$
\[
\Theta_T(rR) D_{\tilde T} = D_{ {\tilde T}^*}(I -  r \tilde R \tilde
T^*)^{-1} (r\tilde R- \tilde T) \quad \mbox{and} \quad I +
\Theta_T(r R) \tilde T^* =  D_{\tilde T^*} (I -r \tilde R\tilde
T^*)^{-1} D_{\tilde T^*}.
\]
\end{lemma}
\noindent \textsf{Proof.} Let $r \in [0, 1)$. Since $\tilde T
D_{\tilde T} = D_{ {\tilde T}^*} \tilde T$ (see equation (3.4) in
Chapter I, \cite{NFBK10}), we have
\[
\begin{split}
\Theta_T(rR) D_{\tilde T} & = \big[- \tilde T +  D_{ {\tilde T}^*}(I
- r \tilde R\tilde T^*)^{-1} r \tilde R D_{\tilde T}\big] D_{\tilde
T}
\\
&= - D_{ {\tilde T}^*}\tilde T+  D_{ {\tilde T}^*}(I -  r \tilde
R\tilde T^*)^{-1} r \tilde R  D_{\tilde T}^2
\\
& = D_{ {\tilde T}^*}(I -  r \tilde R\tilde T^*)^{-1} (-(I -  r
\tilde R\tilde T^*) \tilde T+  r \tilde R D_{\tilde T}^2)
\\
&= D_{ {\tilde T}^*}(I -  r \tilde R\tilde T^*)^{-1} (r\tilde
R-\tilde T).
\end{split}
\]
For the second equality we compute
\[
\begin{split}
D_{\tilde T^*}  (I -r \tilde R\tilde T^*)^{-1}  D_{\tilde T^*} & =
D_{\tilde T^*} (I+(I-r \tilde R\tilde T^*)^{-1}r \tilde R\tilde T^*)
D_{\tilde T^*}
\\
& = D_{\tilde T^*}^2+ D_{\tilde T^*}(I-r \tilde R\tilde T^*)^{-1}r
\tilde R\tilde T^*  D_{\tilde T^*}
\\
& = I - \tilde T \tilde T^*+ D_{\tilde T^*}(I-r \tilde R\tilde
T^*)^{-1} r \tilde R D_{\tilde T} \tilde T^*
\\
& =  I + (- \tilde T + D_{\tilde T^*}(I-r \tilde R\tilde T^*)^{-1} r
\tilde R D_{\tilde T}) \tilde T^*
\\
&=  I  +  \Theta_T(r R) \tilde T^*.
\end{split}
\]
This completes the proof. \qed

We are now ready to prove the main result of this section.

\begin{theorem}\label{factor}
Let $\m H_1$ and $\m H_2$ be two Hilbert spaces and
\[
T =\begin{bmatrix}
A& D_{A^*}L D_B\\
0 & B
\end{bmatrix}:(\op \m H_1) \oplus(\op\m H_2)\to \m H_1 \oplus \m H_2,
\]
be a row contraction on $\m H_1 \oplus \m H_2$ where $A = (A_1,
\ldots, A_n)$ on $\m H_1$ and $B = (B_1, \ldots, B_n)$ on $\m H_2$
are row contractions and $L \in \m B(\m D_B, \m D_{A^*})$ is a
contraction. Then
\[
\Theta_T =(I_\Gamma \otimes \sigma_*^{-1})
 \begin{bmatrix}
\Theta_B & 0 \\
0 & I_{\Gamma \otimes  \m D_{ L^*}}
\end{bmatrix}
(I_\Gamma \otimes J_L)
\begin{bmatrix}
\Theta_A & 0\\ 0& I_{\Gamma \otimes \m D_{ L}}\end{bmatrix}
(I_\Gamma \otimes \sigma),
\]
where $\sigma \in \m B(\m D_T, \m D_A \oplus \m D_L)$ and $\sigma_*
\in \m B(\m D_{T^*}, \m D_{B^*} \oplus \m D_{L^*})$ are unitary
operators as in Theorem \ref{th:NFdefect} and $J_L$ is the
Julia-Halmos matrix corresponding to $L$.
\end{theorem}

\noindent \textsf{Proof.} For each $r \in [0, 1)$, Theorem
\ref{th:NFdefect} and Lemma \ref{le:theta} and the fact that
$(I_\Gamma \otimes \sigma_*)D_{ {\tilde T}^*} = I_{\Gamma} \otimes
\sigma_* D_{T^*}$ yield
\[
(I_{\Gamma} \otimes \sigma_*) \Theta_T(rR) D_{\tilde T} =
\begin{bmatrix}
-\tilde B \tilde L^*D_{\tilde A^*} & D_{\tilde B^*}\\
D_{\tilde L^*}D_{\tilde A^*} & 0
\end{bmatrix}
(I_{\Gamma \otimes \m H} -  r \tilde R\tilde T^*)^{-1} (r\tilde
R-\tilde T).
\]
Now setting $X = D_{A^*}L D_B$, we get
\[
\begin{split}
(I_{\Gamma \otimes \m H}-  r \tilde R\tilde T^*)^{-1} & = \bigg(
\begin{bmatrix}
I_{\Gamma \otimes \m H_1} & 0 \\
0 &       I_{\Gamma \otimes \m H_2}
\end{bmatrix}-
r \begin{bmatrix}
\tilde R & 0 \\
0 &\tilde R
\end{bmatrix}
\begin{bmatrix}
\tilde A^* &  0\\
\tilde X^* & \tilde B^*
\end{bmatrix}\bigg)^{-1}
\\
& = \begin{bmatrix}
I_{\Gamma \otimes \m H_1} -r \tilde R \tilde A^* & 0 \\
- r \tilde R \tilde X^* & I_{\Gamma \otimes \m H_2}- r\tilde R
\tilde B^*
\end{bmatrix}^{-1}
\\
&= \begin{bmatrix}
(I_{\Gamma \otimes \m H_1} -r \tilde R \tilde A^*)^{-1} & 0\\
(I_{\Gamma \otimes \m H_2} - r \tilde R \tilde B^*)^{-1} (r \tilde R
\tilde X^*) (I_{\Gamma \otimes \m H_1} -r \tilde R \tilde A^*)^{-1}&
(I_{\Gamma \otimes \m H_2}- r\tilde R \tilde B^*)^{-1}
\end{bmatrix}
\\
&=\begin{bmatrix}
F& 0 \\
G & H
\end{bmatrix},
\end{split}
\]
where $F =  (I_{\Gamma \otimes \m H_1}- r \tilde R \tilde
A^*)^{-1}$, $G =  H (r \tilde R \tilde X^*)F = (I_{\Gamma \otimes \m
H_2}- r\tilde R \tilde B^*)^{-1} (r \tilde R \tilde X^*) (I_{\Gamma
\otimes \m H_1}-r \tilde R \tilde A^*)^{-1}$ and $H = (I_{\Gamma
\otimes \m H_2}- r \tilde R \tilde B^*)^{-1}$. Therefore
\[
\begin{split}
(I_{\Gamma} \otimes \sigma_*) \Theta_T(rR) D_{\tilde T} & =
\begin{bmatrix}
-\tilde B \tilde L^*D_{\tilde A^*} & D_{\tilde B^*}\\
D_{\tilde L^*}D_{\tilde A^*} & 0
\end{bmatrix}
\begin{bmatrix}
F& 0 \\
G &H
\end{bmatrix}
\begin{bmatrix}
r \tilde R - \tilde A & - \tilde X \\
0 &r \tilde R - \tilde B \\
\end{bmatrix}
\\
& =
\begin{bmatrix}
-\tilde B \tilde L^*D_{\tilde A^*}&   D_{\tilde B^*}\\
D_{\tilde L^*}D_{\tilde A^*} & 0
\end{bmatrix}
\begin{bmatrix}
F (r \tilde R - \tilde A)& - F \tilde X\\
G (r \tilde R - \tilde A) & -G \tilde X + H(r \tilde R - \tilde B)\\
\end{bmatrix}
\\
& =
\begin{bmatrix}
C_{11}(r) & C_{12}(r)\\
C_{21}(r) & C_{22}(r)
\end{bmatrix} \in \m B((\Gamma \otimes \m H _1) \oplus (\Gamma
\otimes \m H_2)),
\end{split}
\]
where $C_{11}(r) = -\tilde B \tilde L^*D_{\tilde A^*} F  (r \tilde R
- \tilde A)+  D_{\tilde B^*} G (r \tilde R - \tilde A)$, $C_{12}(r)
= \tilde B \tilde L^*D_{\tilde A^*} F   \tilde X- D_{\tilde
B^*}G\tilde X + D_{\tilde B^*}H(r \tilde R - \tilde B)$, $C_{21}(r)
= D_{\tilde L^*}D_{\tilde A^*}  F (r \tilde R - \tilde A)$ and
$C_{22}(r) = - D_{\tilde L^*}D_{\tilde A^*}  F \tilde X$. Further,
we compute
\[
\begin{split}
C_{11}(r) & = -\tilde B \tilde L^*D_{\tilde A^*} F  (r \tilde R -
\tilde A)+ D_{\tilde B^*} (H r \tilde R \tilde X^* F) (r \tilde R -
\tilde A)
\\
&= -\tilde B \tilde L^*D_{\tilde A^*} F  (r \tilde R - \tilde A)+
D_{\tilde B^*} H  (r \tilde R \tilde X^*) F (r \tilde R- \tilde
A)
\\
&= -\tilde B \tilde L^*D_{\tilde A^*} F  (r \tilde R- \tilde A)+
D_{\tilde B^*} H  (r \tilde R D_{\tilde B} \tilde L^* D_{\tilde
A^*}) F (r \tilde R- \tilde A)
\\
&= [ -\tilde B+  D_{\tilde B^*} H (r \tilde R D_{\tilde B})] \tilde
L^*D_{\tilde A^*} F (r \tilde R- \tilde A)
\\
&= [ -\tilde B+  D_{\tilde B^*} (I_{\Gamma \otimes \m H_2}- r\tilde
R\tilde B^*)^{-1} r \tilde R D_{\tilde B}] \tilde L^*D_{\tilde A^*}
F (r \tilde R- \tilde A)
\\
&= \Theta_B(r R)\tilde L^*D_{\tilde A^*} F  (r \tilde R- \tilde A)
\\
&= \Theta_B(r R)\tilde L^*D_{\tilde A^*}  (I_{\Gamma \otimes \m
H_1}-r \tilde R\tilde A^*)^{-1}  (r \tilde R- \tilde A)
\\
&= \Theta_B(r R)\tilde L^* \Theta_A(r R) D_{\tilde A},
\end{split}
\]
where the last equality follows from Lemma \ref{le:theta}. Also
\[
\begin{split}
C_{12}(r)&= \tilde B \tilde L^*D_{\tilde A^*} F   \tilde X-
D_{\tilde B^*}G\tilde X + D_{\tilde B^*}H(r \tilde R - \tilde B)
\\
&= \tilde B \tilde L^*D_{\tilde A^*} F   \tilde X-  D_{\tilde B^*} H
(r \tilde R \tilde X^*) F  \tilde X + D_{\tilde B^*}H(r \tilde R -
\tilde B)
\\
& = \tilde B \tilde L^*D_{\tilde A^*} F   \tilde X-  D_{\tilde B^*}
H (r \tilde R D_{\tilde B} \tilde L^*D_{\tilde A^*} ) F  \tilde X +
D_{\tilde B^*}H(r \tilde R - \tilde B)
\\
& = -[-\tilde B +  D_{\tilde B^*} H  r \tilde R D_{\tilde B} ]
\tilde L^*D_{\tilde A^*} F   \tilde X + D_{\tilde B^*}H(r \tilde R -
\tilde B)
\\
& = -[-\tilde B +  D_{\tilde B^*} (I_{\Gamma \otimes \m H_2}-
r\tilde R\tilde B^*)^{-1} r \tilde R D_{\tilde B} ] \tilde
L^*D_{\tilde A^*} F \tilde X + D_{\tilde B^*}H(r \tilde R - \tilde
B)
\\
& = - \Theta_B(r R) \tilde L^*D_{\tilde A^*} F   D_{\tilde
A^*}\tilde L D_{\tilde B} + D_{\tilde B^*} (I_{\Gamma \otimes \m
H_2}- r\tilde R\tilde B^*)^{-1}  (r \tilde R - \tilde B)
\\
&= - \Theta_B(r R) \tilde L^*D_{\tilde A^*} F   D_{\tilde A^*}\tilde
L D_{\tilde B} + \Theta_B(r R)D_{\tilde B} \quad \quad (\mbox{by
Lemma }~\ref{le:theta})
\\
&= \Theta_B(r R)[ - \tilde L^* (D_{\tilde A^*} F D_{\tilde
A^*})\tilde L  +I_{\m D_{\tilde B}}]D_{\tilde B}
\\
&= \Theta_B(r R)[ - \tilde L^* (D_{\tilde A^*}  (I_{\Gamma \otimes
\m H_1} - r \tilde R\tilde A^*)^{-1}  D_{\tilde A^*}) \tilde L
+I_{\m D_{\tilde B}}]D_{\tilde B}
\\
& = \Theta_B(r R)[ - \tilde L^* (I_{\Gamma \otimes \m H_1}  +
\Theta_A(r R) \tilde A^*)\tilde L  +I_{\m D_{\tilde B}}]D_{\tilde B}
\quad \quad (\mbox{by Lemma }~\ref{le:theta})
\\
&=  \Theta_B(r R) (- \tilde L^* \Theta_A(r R) \tilde A^*\tilde L +
D_{\tilde L}^2)D_{\tilde B},
\end{split}
\]
and
\[
\begin{split}
C_{21}(r) & = D_{\tilde L^*}D_{\tilde A^*} F (r \tilde R- \tilde A)
= D_{\tilde L^*}D_{\tilde A^*}(I_{\Gamma \otimes \m H_1}-r \tilde
R\tilde A^*)^{-1} (r \tilde R- \tilde A) = D_{\tilde L^*} \Theta_A(r
R) D_{\tilde A},
\end{split}
\]
and finally
\[
\begin{split}
C_{22}(r) & = - D_{\tilde L^*}D_{\tilde A^*}F \tilde X = - D_{\tilde
L^*} D_{\tilde A^*} (I_{\Gamma \otimes \m H_1}-r \tilde R\tilde
A^*)^{-1} D_{\tilde A^*}  \tilde L D_{\tilde B}
\\
&= - D_{\tilde L^*} [I_{\Gamma \otimes \m H_1}  +  \Theta_A(r R)
\tilde A^* ] \tilde L D_{\tilde B} \quad \quad (\mbox{by Lemma}~\ref{le:theta})
\\
& = - D_{\tilde L^*} \Theta_A(r R) \tilde A^* \tilde L D_{\tilde B} - D_{\tilde L^*} \tilde L D_{\tilde B} \\
& = - D_{\tilde L^*} \Theta_A(r R) \tilde A^* \tilde L D_{\tilde B}
- \tilde L D_{\tilde L} D_{\tilde B}.
\end{split}
\]
This implies that
\[
\begin{split}
(I_\Gamma \otimes \sigma_*)\Theta_T(rR) D_{\tilde T} & =
\begin{bmatrix}
\Theta_B(r R)\tilde L^* \Theta_A(r R) D_{\tilde A} &
\Theta_B(r R)(- \tilde L^* \Theta_A(r R) \tilde A^*\tilde L + D_{\tilde L}^2) D_{\tilde B}\\
 D_{\tilde L^*}  \Theta_A(r R) D_{\tilde A}&
 - D_{\tilde L^*} \Theta_A(r R) \tilde A^* \tilde L D_{\tilde B} - \tilde L D_{\tilde L} D_{\tilde B}\\
\end{bmatrix}
\\
& = \begin{bmatrix} \Theta_B(r R) \tilde L^* \Theta_A(rR) &
\Theta_B( r R) D_{\tilde L}
\\
D_{\tilde L^*}  \Theta_A( r R) & - \tilde L \\
\end{bmatrix}
\begin{bmatrix}
D_{\tilde A} &
- \tilde A^*\tilde L D_{\tilde B} \\
0 &
D_{\tilde L} D_{\tilde B}\\
\end{bmatrix}
\\
&= \begin{bmatrix} \Theta_B(r R) \tilde L^* \Theta_A(rR) &
\Theta_B( r R) D_{\tilde L} \\
D_{\tilde L^*}  \Theta_A( r R) & - \tilde L
\end{bmatrix} (I_\Gamma \otimes \sigma) D_{\tilde T} \quad ({\rm
by ~ Theorem ~\ref{th:NFdefect}}),
\end{split}
\]
and we conclude that
\[
(I_\Gamma \otimes \sigma_*)\Theta_T(rR) =
\begin{bmatrix}
\Theta_B(r R) \tilde L^* \Theta_A(rR) &
\Theta_B( r R) D_{\tilde L} \\
D_{\tilde L^*}  \Theta_A( r R) & - \tilde L
\end{bmatrix}
(I_\Gamma \otimes \sigma).
\]
We may rewrite this as
\[
(I_\Gamma \otimes \sigma_*)\Theta_T(rR) (I_\Gamma \otimes
\sigma^{-1}) =
\begin{bmatrix}
\Theta_B(r R) \tilde L^* \Theta_A(rR) &
\Theta_B( r R) D_{\tilde L} \\
D_{\tilde L^*}  \Theta_A( r R) & - \tilde L
\end{bmatrix}.
\]
Finally, we observe that
\[
\begin{split}
\begin{bmatrix}
\Theta_B(rR)\tilde L^* \Theta_A( rR)  &  \Theta_B(r R) D_{\tilde L} \\
D_{\tilde L^*}  \Theta_A( rR)  & - \tilde L
\end{bmatrix}
&=
\begin{bmatrix}
\Theta_B(rR) &  0 \\
0 & I_{\m D_{\tilde L^*}}
\end{bmatrix}
\begin{bmatrix}
 \tilde L^* &D_{\tilde L} \\
 D_{\tilde L^*}  & - \tilde L \\
\end{bmatrix}
 \begin{bmatrix}
\Theta_A( r R) &  0 \\
0   & I_{\m D_{\tilde L}}
\end{bmatrix}
\\
&=
\begin{bmatrix}
\Theta_B(r R)  &   0\\
0 & I_{\Gamma \otimes  \m D_{ L^*}}\\
\end{bmatrix}
\big( I_\Gamma \otimes J_L)
\begin{bmatrix}
\Theta_A(rR)  &  0 \\
 0  & I_{\Gamma \otimes \m D_{ L}} \\
\end{bmatrix},
\end{split}
\]
so that the resulting formula is
\[
\begin{split} \Theta_T(rR)=(I_\Gamma \otimes \sigma_*^{-1})
 \begin{bmatrix}
\Theta_B(r R)&   0 \\
  0 & I_{\Gamma \otimes  \m D_{ L^*}} \\
\end{bmatrix}
(  I_\Gamma \otimes J_L)
\begin{bmatrix}
\Theta_A( rR)  &   0 \\
0  & I_{\Gamma \otimes \m D_{ L}}\\
\end{bmatrix}
(I_\Gamma \otimes \sigma).
\end{split}
\]
The result now follows by passing to the strong operator topology
limit as $r \to 1$. \qed

In the following, we prove that the Julia-Halmos matrix factor $J_L$
in the factorization of the above theorem is canonical. The proof is
similar to the one for $n = 1$ case by Sz.-Nagy and Foias (see
Theorem 3, page 209-212, \cite{NF67}). We only sketch the main ideas
and refer to \cite{NF67} for full proof details.

\begin{theorem}\label{converse.factor}
Let $\m H_1, \m H_2, \m F$ and $\m F_*$ be Hilbert spaces and $A =
(A_1, \ldots, A_n)$ and $ B = (B_1, \ldots, B_n) $ be row
contractions on $\m H_1$ and $\m H_2$, respectively. Let $w \in \m
B(\m D_{A^*} \oplus \m F, \m D_B \oplus \m F_*)$ be a unitary
operator and
\[
\Theta =  \begin{bmatrix} \Theta_B  &  0
\\
0 &
I_{\Gamma \otimes  \m F_*} \\
\end{bmatrix}
(I_\Gamma\otimes  w )
\begin{bmatrix}
\Theta_A  &  0  \\
0 & I_{\Gamma \otimes  \m F} \\
\end{bmatrix}
: \Gamma \otimes (\m D_A \oplus \m F) \to \Gamma \otimes (\m
D_{B^*}\oplus \m F_*)
\]
be a purely contractive multi-analytic operator. Then $\Theta$ and
$\Theta_T$ coincide where
\[
T =\begin{bmatrix}
A& D_{A^*} (P_{\m D_{A^*}} w^*|_{\m D_B})D_B\\
0 & B\\
\end{bmatrix}:(\op \m H_1) \oplus(\op\m H_2)\to \m H_1 \oplus \m
H_2.
\]
\end{theorem}

\noindent \textsf{Proof.} Let $w^* = \begin{bmatrix}
L& M\\
N& K\\
\end{bmatrix}$ where $L = P_{\m D_{A^*}} w^*|_{\m D_B} \in \m B(\m D_B, \m
D_{A^*})$, $M \in \m B(\m F_*, \m D_{A^*})$, $N \in \m B(\m D_B, \m
F)$ and $K \in \m B(\m F_*, \m F)$ are contractions. Define $ \m
F\p:= \m F\ominus \ov{N\m D_B} $ and $ \m F_*\p:= \m F_*\ominus
\ov{M^*\m D_{A^*}} $. Following the same line of argument as in the
proof of the first part of Theorem 3 in \cite{NF67} we have
\[
w \m F\p  = \m F_*\p.
\]
In particular, for each $ f\p \in \m F\p (\subset \m F)$ we have $w
f\p \in \m F_*\p (\subset\m F_*)$ and
\begin{align*}
\Theta (e_\emptyset \otimes f\p)
&=
\begin{bmatrix}
\Theta_B  &  0 \\
0 & I_{\Gamma \otimes  \m F_*} \\
\end{bmatrix}
(I_\Gamma\otimes  w )
\begin{bmatrix}
\Theta_A &  0\\
0 & I_{\Gamma \otimes  \m F} \\
\end{bmatrix}
\big (e_\emptyset \otimes (0\oplus f\p)\big ) \\
&=
\begin{bmatrix}
\Theta_B &  0\\
0 & I_{\Gamma \otimes  \m F_*}
\end{bmatrix}
(e_\emptyset\otimes  (0 \oplus w f\p))\\
&= 0 \oplus (e_\emptyset \otimes  w f\p).
\end{align*}
Then $ \|P_{e_\emptyset\otimes (\m D_{B^*} \oplus \m F_*)}\Theta
(e_\emptyset \otimes f\p)\|^2 =\| e_\emptyset \otimes  w f\p\|^2 =
\|f\p\|^2 $. Since $ \Theta $ is a purely contractive, $f\p = 0$,
that is, $\m F\p = \{0\}$ and hence $\m F_*\p = \{0\}$. Hence $
\ov{N\m D_B} = \m F$ and $\ov{M^*\m D_{A^*}} = \m F_*$.
Consequently, $U \in \m B(\m F, \m D_L)$ and $V \in \m B(\m F_*, \m
D_{L^*})$ defined by
\[
U(N x) = D_L x \quad \mbox{and} \quad V (M^*y)= D_{L^*}y \quad \quad
(x \in \m D_B, y \in \m D_{A^*}),
\]
are unitary operators. Also
\[
N^* = D_L|_{\m D_L} U.
\]
Then
\[ w=
\begin{bmatrix}
L^*&   N^*\\
M^*& K^*
\end{bmatrix} = \begin{bmatrix}
L^*&  D_L|_{\m D_L} U\\
M^*& V^*K_1 U
\end{bmatrix}
 = v^* J u,
\]
where $K_1 = VK^*U^* \in \m B(\m D_L, \m D_{L^*})$, $u =
\begin{bmatrix}
I_{\m D_{A^*}}& 0\\
0 & U\end{bmatrix}$, $v = \begin{bmatrix}
I_{\m D_{B}}& 0\\
0 & V\end{bmatrix}$ and $J = \begin{bmatrix} L^*&  D_L|_{\m D_L}
\\ D_{L^*} & K_1
\end{bmatrix}$.

\noindent Since $J \in \m B(\m D_{A^*} \oplus \m D_L, \m D_B \oplus
\m D_{L^*})$ is a unitary operator we have (see page 211 in
\cite{NF67}) $ K_1 = -L|_{\m D_L} $. Now for $ u\p :=
\begin{bmatrix} I_{\m D_{A}}& 0\\ 0 & U\end{bmatrix}$ and $v\p :=
\begin{bmatrix} I_{\m D_{B^*}}& 0\\ 0 & V\end{bmatrix}$, we have
\[
(I_\Gamma  \otimes v\p )
\begin{bmatrix}
  \Theta_B  &  0 \\
0 & I_{\Gamma \otimes  \m F_*} \\
\end{bmatrix}
=
\begin{bmatrix}
\Theta_B &   0 \\
0 & I_{\Gamma \otimes  \m D_{L^*}} \\
\end{bmatrix}
(I_{\Gamma}\otimes v ),
\] and
\[
(I_\Gamma  \otimes u )
\begin{bmatrix}
\Theta_A  &  0\\
0 & I_{\Gamma \otimes  \m F} \\
\end{bmatrix}
=
\begin{bmatrix}
\Theta_A  &   0 \\
0 & I_{\Gamma \otimes  \m D_{L}} \\
\end{bmatrix}
 (I_{\Gamma}\otimes u\p ).
\]
This implies that

\[
\begin{split}
(I_\Gamma  \otimes v\p ) \Theta(I_{\Gamma}\otimes {u\p}^* )
&=
(I_\Gamma  \otimes v\p )
\begin{bmatrix}
\Theta_B  &  0 \\
0 & I_{\Gamma \otimes  \m F_*} \\
\end{bmatrix}
(I_{\Gamma}\otimes w )
\begin{bmatrix}
\Theta_A  &  0  \\
 0 & I_{\Gamma \otimes  \m F} \\
 \end{bmatrix}
(I_{\Gamma}\otimes {u\p}^* )
\\
&=
\begin{bmatrix}
\Theta_B &   0 \\
 0 & I_{\Gamma \otimes  \m D_{L^*}} \\
\end{bmatrix}
(I_\Gamma  \otimes v)(I_{\Gamma}\otimes w )(I_{\Gamma}\otimes u^* )
\begin{bmatrix}
\Theta_A &   0 \\
0 & I_{\Gamma \otimes  \m D_{L}} \\
\end{bmatrix}
\\
&=
\begin{bmatrix}
\Theta_B &   0 \\
0 & I_{\Gamma \otimes  \m D_{L^*}} \\
\end{bmatrix}
(I_{\Gamma}\otimes J)
\begin{bmatrix}
\Theta_A  &   0 \\
 0 & I_{\Gamma \otimes  \m D_{L}} \\
\end{bmatrix}.
\end{split}
\]
Since $L = P_{\m D_{A^*}} w^*|_{\m D_B}$ is a contraction, Theorem
\ref{th:NFblock} shows that the $n$-tuple $T$ defined by $T
=\begin{bmatrix}
A& D_{A^*}L D_B\\
0 & B\end{bmatrix}$ is a row contraction on $ \m H_1 \oplus \m H_2$
and Theorem \ref{factor} implies that
\[
(I_{\Gamma}\otimes  \sigma_*)\Theta_T (I_{\Gamma} \otimes
\sigma^{-1}) =  \begin{bmatrix}
\Theta_B &   0 \\
0 &
I_{\Gamma \otimes  \m D_{L^*}} \\
\end{bmatrix}
(I_{\Gamma}\otimes J_L )
\begin{bmatrix}
\Theta_A &0 \\
0 & I_{\Gamma \otimes  \m D_{L}} \\
\end{bmatrix}
\]
where $ \sigma: \m D_T \to \m D_A \oplus \m D_L$ and $\sigma_*: \m
D_{T^*} \to \m D_{B^*} \oplus \m D_{L^*}$ are unitary operators as
in Theorem \ref{th:NFblock}. Therefore, $(I_\Gamma \otimes v\p )
\Theta(I_{\Gamma}\otimes {u\p}^* )= (I_{\Gamma}\otimes
\sigma_*)\Theta_T (I_{\Gamma} \otimes \sigma^{-1})$, that is,
$\Theta_T$ coincides with $ \Theta $. \qed

\section{Factorizations of Characteristic Functions of Constrained Row Contractions}

The main objective of this section is to study factorizations of
characteristic functions of row contractions in noncommutative
varieties. The notion of a noncommutative variety was introduced by
G. Popescu in \cite{Po06a}.

We first recollect some basic definitions, notations, and results
that will be used subsequently. For details, we refer to
\cite{Po06a}, \cite{P10} and the references therein. Let $ \m P_J
\subset \m F_n^\infty$ be a family of noncommutative polynomials and
$J$ be the WOT-closed two sided ideal of $ \m F_n^\infty$ generated
by $\m P_J$. In what follows, we always assume that $J \neq \m
F_n^\infty$. Then
\[
\m M_J := \ov{\mbox{span}} \{\phi \otimes \psi: \phi \in J, \psi\in
\Gamma\} \quad \mbox{and} \quad \m N_J : = \Gamma\ominus \m M_J,
\]
are proper joint $(L_1, \ldots, L_n)$ and $(L_1^*, \ldots, L_n^*)$
invariant subspaces of $\Gamma$, respectively. Define constrained
left creation operators and constrained right creation operators on
$\m N_J$ by
\[
V_j := P_{\m N_J }L_j|_{P_{\m N_J }}\quad \quad {\rm and} \quad
\quad W_j := P_{\m N_J }R_j|_{P_{\m N_J }} \quad \quad (j =1,
\ldots, n),
\]
respectively.

Let $\m E$ and $\m E_*$ be Hilbert spaces and $M\in \z(\m N_J
\otimes \m E, \m N_J \otimes \m E_* )$. Then $M$ is said to be
\textit{ constrained multi-analytic operator} if
\[
M(V_j \otimes I_{\m E})= (V_j \otimes I_{\m E_*})M \quad \quad \quad
(j =1,\ldots, n).
\]
We say that $M\in \z(\m N_J \otimes \m E, \m N_J \otimes \m E_* ) $
is \textit{purely contractive} if $M$ is a contraction and
$e_\emptyset \in \m N_J$ and
\[
\|P_{e_\emptyset \otimes \m E_*} M(e_\emptyset \otimes \eta)\| <
\|\eta \| \quad \quad (\eta \neq 0,  \eta \in \m E).
\]
It has been shown by Popescu \cite{Po06a} that the set of all
constrained multi-analytic operators in $\z(\m N_J \otimes \m E, \m
N_J \otimes \m E_* )$ coincides with
\[
\m W(W_1,\ldots, W_n)~ \bar\otimes ~\z( \m E,\m E_* ) = P_{\m N_J
\otimes \m E_*}[ \m R_n^\infty~\bar \otimes~\z( \m E,\m E_* )
]|_{P_{\m N_J\otimes \m E}}.
\]
where $\m W(W_1,\ldots, W_n)$ is the WOT-closed algebra generated by
$\{I, W_1,\ldots, W_n\}$ and $\m R_n^\infty = U^* \m F_n^\infty U$
and $U$ is the flipping operator.

A row contraction $T= (T_1, \ldots, T_n)$ on $\m H$ is said to be
\textit{$J$-constrained row contraction}, or simply
\textit{constrained row contraction} if $J$ is clear from the
context, if
\[
p(T_1, \ldots, T_n) = 0 \quad \quad \quad (p \in \m P _J).
\]

The \textit{constrained characteristic function} (see Popescu
\cite{Po06a}) of a constrained row contraction $T= (T_1, \ldots,
T_n)$ on a Hilbert space $\m H$ is the constrained multi-analytic
operator $\Theta_{J, T}:\m N_J\otimes\m D_{T} \to \m N_J\otimes\m
D_{T^*}$ defined by
\[
\Theta_{J, T}= P_{\m N_J \otimes \m D_{T^*}}\Theta_T|_{\m N_J
\otimes \m D_{T}}.
\]
Since $ \m N_J $ is a joint $(R_1^* \otimes I_{\m D_{T^*}}, \ldots,
R_n^* \otimes I_{\m D_{T^*}})$ invariant subspace and $W_j = P_{\m
N_J }(R_j \otimes I_{\m D_{T^*}})|_{P_{\m N_J }},~~ j =1, \ldots,
n$, it follows that (see \cite{Po06a})
\begin{gather}
\Theta_T^*(\m N_J \otimes \m D_{T^*}) \subset  \m N_J \otimes \m
D_{T} \quad \mbox{and} \quad \Theta_T(\m M_J \otimes \m D_{T})
\subset  \m M_J \otimes \m D_{T^*}.
\end{gather}

\textsf{From here onwards to maintain simplicity of notations, we
often omit the subscript $J$.}

Now we are ready to prove a factorization of constrained
characteristic functions corresponding to upper triangular
constrained row contractions.

\begin{theorem}\label{cons.factor.}
Let $T =
\begin{bmatrix}
A&  D_{A^*}L D_B\\
0 & B
\end{bmatrix}$ be a constrained row contraction on $\m H_1 \oplus \m
H_2$ where $A = (A_1, \ldots, A_n)$ on $\m H_1$ and $B = (B_1,
\ldots, B_n)$ on $\m H_2$ are row contractions and $L \in \m B(\m
D_B, \m D_{A^*})$ is a contraction. Then $A$ and $B$ are also
constrained row contractions and
\begin{align*}
\Theta_{J, T}= (I_{\m N}  \otimes \sigma_*^{-1})
\begin{bmatrix}
\Theta_{J, B} &  0 \\
 0 &I_{\m N \otimes  \m D_{ L^*}} \\
\end{bmatrix}
(I_{\m N}\otimes J_L)
\begin{bmatrix}
\Theta_{J,A}  & 0 \\
0  &  I_{\m N \otimes \m D_{ L}} \\
\end{bmatrix}
(I_{\m N} \otimes \sigma)
\end{align*}
where $\sigma \in \m B(\m D_T, \m D_A \oplus \m D_L)$ and $\sigma_*
\in \m B(\m D_{T^*}, \m D_{B^*} \oplus \m D_{L^*})$ are unitary
operators as in Theorem \ref{th:NFblock}.
\end{theorem}
\noindent \textsf{Proof.} It is straightforward to verify that $A$
and $B$ are constrained row contractions. For the remaining part,
first we observe that
\[
\Theta_{J, T} = P_{\m N \otimes \m D_{T^*}} (I_\Gamma  \otimes
\sigma_*^{-1})
\begin{bmatrix}
\Theta_B & 0\\
0 & I_{\Gamma \otimes  \m D_{ L^*}}
\end{bmatrix}
(I_\Gamma \otimes  J_L )
\begin{bmatrix}
\Theta_A &  0 \\
0 & I_{\Gamma \otimes \m D_{ L}}
\end{bmatrix}
(I_\Gamma \otimes \sigma)|_{\m N \otimes \m D_T}.
\]
Since $P_{\m N \otimes \m D_{T^*}} (I_\Gamma  \otimes \sigma_*^{-1})
= (I_{\m N} \otimes  \sigma_*^{-1}) P_{\m N \otimes (\m D_{B^*}
\oplus \m D_{L^*})}$ and
\[
P_{\m N \otimes (\m D_{B^*}
\oplus \m D_{L^*})}
\begin{bmatrix}
\Theta_B & 0\\
0 & I_{\Gamma \otimes  \m D_{ L^*}}
\end{bmatrix}
= P_{\m N \otimes (\m D_{B^*} \oplus \m D_{L^*})}
\begin{bmatrix}
\Theta_B & 0\\
0 & I_{\m N \otimes  \m D_{ L^*}}
\end{bmatrix}
P_{\m N \otimes (\m D_{B} \oplus \m D_{L^*})},
\]
and $P_{\m N \otimes (\m D_{B} \oplus \m D_{L^*})} (I_{\Gamma}
\otimes J_L) = I_{\m N} \otimes J_L = (I_{\m N} \otimes J_L) P_{\m N
\otimes (\m D_{A^*} \oplus \m D_{L})}$, and
\[
P_{\m N \otimes (\m D_{A^*} \oplus \m D_{L})}
\begin{bmatrix}
\Theta_A &  0 \\
0 & I_{\Gamma \otimes \m D_{ L}}
\end{bmatrix}
=
P_{\m N \otimes (\m D_{A^*} \oplus \m D_{L})}
\begin{bmatrix}
\Theta_A &  0 \\
0 & I_{\m N \otimes \m D_{ L}}
\end{bmatrix}
P_{\m N \otimes (\m D_{A} \oplus \m D_{L})},
\]
we have the required equality. \qed

We now state a similar result to Theorem \ref{converse.factor} for
constrained row contractions. We omit the proof, which uses similar
techniques to the proof of Theorem \ref{converse.factor} (and
Theorem 3 in \cite{NF67}).

\begin{theorem}\label{cons.factor.-conv}
Let $\m H_1, \m H_2, \m F$ and $\m F_*$ be Hilbert spaces and $A =
(A_1, \ldots, A_n)$ and $ B = (B_1, \ldots, B_n) $ be constrained
row contractions on $\m H_1$ and $\m H_2$, respectively, and
$e_\emptyset \in \m N$. Let $w \in \m B(\m D_{A^*} \oplus \m F, \m
D_B \oplus \m F_*)$ be a unitary operator and $T =\begin{bmatrix}
A& D_{A^*} (P_{\m D_{A^*}} w^*|_{\m D_B})D_B\\
0 & B\\
\end{bmatrix}$ be a constrained row contraction and
\[
\Theta =  \begin{bmatrix}
  \Theta_{J,B}  &   0\\
0 & I_{\m N \otimes  \m F_*} \\
\end{bmatrix}
(I_{\m N}\otimes  w )
\begin{bmatrix}
 \Theta_{J,A}  &  0 \\
0 & I_{\m N \otimes  \m F} \\
\end{bmatrix}
\]
be a purely contractive constrained multi-analytic operator. Then
$\Theta$ coincides with $\Theta_{J,T}$.
\end{theorem}

A particularly important example of noncommutative variety is the
one given by $\m P_{J_c} = \{L_i L_j - L_j L_i: i, j = 1, \ldots,
n\}$. In this case $\m N_{J_c} = \Gamma_s$ is the symmetric Fock
space, $V_j =  P_{\Gamma_s}L_j|\Gamma_s, j=1, \ldots, n$, are the
creation operators on $\Gamma_s$ (see \cite{B"}, \cite{Po06a}).
Moreover, one can identify $(V_1, \ldots, V_n)$ on $\Gamma_s$ with
the multiplication operator tuple $(M_{z_1}, \ldots, M_{z_n})$ on
the Drury-Arveson space $H^2_n$. Recall that the Drury-Arveson space
is a reproducing kernel Hilbert space with kernel function $k:
\mathbb B^n \times \mathbb B^n \to \mathbb C$ defined by
\[
k(z, w) = (1- \langle z,w \rangle_{\mathbb C^n})^{-1} \quad \quad
(z, w \in \mathbb B^n).
\]
Under this identification, the set of constrained multi-analytic
operators $P_{\Gamma_s}\m F_n^\infty|_{\Gamma_s}$ corresponds to the
multiplier algebra of $H^2_n$.

\noindent  Note too that a row contraction $T= (T_1, \ldots, T_n)$
on $\m H$ is a constrained row contraction if and only if $T$ is a
commuting row contraction, that is, $T_i T_j = T_j T_i$, $i, j = 1,
\ldots, n$. In this case, we can identify the constrained
characteristic function $\Theta_{J_c, T}= P_{\m N_{J_c} \otimes \m
D_{T^*}}\Theta_T|_{\m N_{J_c} \otimes \m D_{T}}$ with the bounded
operator-valued analytic function $\theta_T$ on $\mathbb B^n$
defined by (see \cite{Po06a}, \cite{B"} and \cite{BES})
\[
\theta_{T}(z) = -T + D_{T^*} (I_{\m H}- \sum_{i=1}^n z_i
T_i^*)^{-1}ZD_T \quad \quad (z \in \mathbb B^n),
\]
where $Z = (z_1 I_{\m H}, \ldots, z_n I_{\m H})$, $z \in \mathbb
B^n$.

In this setting, Theorems \ref{cons.factor.} and
\ref{cons.factor.-conv} can be stated as follows:

\begin{theorem}
Let $T =\begin{bmatrix}
A& D_{A^*} L D_B\\
0 & B
\end{bmatrix}$ be a commuting row contraction on $ \m H_1 \oplus \m
H_2$ where $A$ and $B$ are commuting row contractions on $\m H_1$
and $\m H_2$, respectively, and $L \in \m B(\m D_B, \m D_{A^*})$ is
a contraction. Then $\theta_T$ coincides with
\[
\begin{bmatrix}
  \theta_{B}  &   0\\
0 & I_{H^2_n \otimes  \m D_{L^*}} \\
\end{bmatrix}
(I_{H^2_n}\otimes  J_L )
\begin{bmatrix}
 \theta_{A}  &  0 \\
0 & I_{H^2_n \otimes  \m D_L} \\
\end{bmatrix}.
\]
Moreover, if $\hat{T} =\begin{bmatrix}
A& D_{A^*} (P_{\m D_{A^*}} w^*|_{\m D_B})D_B\\
0 & B
\end{bmatrix}$ is a commuting row contraction for some unitary
operator $w \in \m B(\m D_{A^*} \oplus \m F, \m D_B \oplus \m F_*)$
and Hilbert spaces $\m F$ and $\m F_*$, and if
\[
\theta =  \begin{bmatrix}
  \theta_{B}  &   0\\
0 & I_{H^2_n \otimes  \m F_*} \\
\end{bmatrix}
(I_{H^2_n} \otimes  w )
\begin{bmatrix}
 \theta_{A}  &  0 \\
0 & I_{H^2_n \otimes  \m F} \\
\end{bmatrix}
\]
is a purely contractive multiplier then $\theta$ coincides with
$\theta_{\hat T}$.
\end{theorem}

Now let $\m H_1$ be a closed subspace of a Hilbert space $\m H$ and
$T = (T_1, \ldots, T_n)$ be an $n$-tuple on $\m H$. Let $\m H_1$ be
a joint $T$ invariant subspace of $\m H$ (that is, $T_i \m H_1
\subseteq \m H_1$ for all $i = 1, \ldots, n$) and $\m H_2 = \m H
\ominus  \m H_1$. Then we can represent, with respect to $\m H = \m
H_1 \oplus \m H_2$, $T_j$ as an upper triangular operator matrix
\[
T_j =\begin{bmatrix}
A_j&  X_j\\
0 & B_j
\end{bmatrix},
\]
where $A_j = T_j|_{\m H_1} \in \z(\m H_1)$, $B_j = P_{\m H_2}
T_j|_{\m H_2} \in \z (\m H_2)$ and $X_j = P_{\m H_1} T_j|_{\m H_2}
\in \z(\m H_2, \m H_1)$, $j = 1, \ldots, n$. In other words
\begin{equation}\label{eq:2b2T}
T =\begin{bmatrix}
A&  X\\
0 & B\\
\end{bmatrix}:(\op \m H_1) \oplus(\op \m H_2)\to \m H_1 \oplus \m
H_2,
\end{equation}
where $A = (A_1, \ldots, A_n)\in \z (\op \m H_1, \m H_1)$, $B =
(B_1, \ldots, B_n) \in \z( \op \m H_2, \m H_2)$ and $X = (X_1,
\ldots, X_n) \in \z(\op \m H_2,\m H_1 )$.

Conversely, let $T$ be a row operator on $\m H$ and $\m H_1$ and $\m
H_2$ be closed subspaces of $\m H$. If $T$ admits an upper
triangular representation as in (\ref{eq:2b2T}) for some row
operators $A = (A_1, \ldots, A_n)\in \z (\op \m H_1, \m H_1)$, $B =
(B_1, \ldots, B_n) \in \z( \op \m H_2, \m H_2)$ and $X = (X_1,
\ldots, X_n) \in \z(\op \m H_2,\m H_1 )$ then $\m H_1$ is a joint
$T$-invariant subspace of $\m H$. In other words, $T$ has a
non-trivial joint invariant subspace if and only if $T$ admits an
upper triangular representation as in (\ref{eq:2b2T}). This is also
equivalent to the regular factorizations of the characteristic
function $\Theta_T$ in terms of $\Theta_A$ and $\Theta_B$ (see
Sz.-Nagy and Foias \cite{NFBK10} for $n =1$ case and Popescu
\cite{Po06b} for general case). It is not known, in the general
case, how one relates regular factorizations of characteristic
functions and the one obtained in this paper. We do not know the
answer even if $n = 1$.

\vspace{0.2in}

\noindent\textit{Acknowledgement:} The first author's research
work is supported by DST-INSPIRE Faculty Fellowship
No. DST/INSPIRE/04/2014/002624. The second author's research work
is supported by NBHM Post Doctoral Fellowship
No. 2/40(50)/2015/ R \& D  - II/11569. The third author is supported in
part by NBHM (National Board of Higher Mathematics, India) grant
NBHM/R.P.64/2014.


\end{document}